\documentclass[12pt]{amsart}

\usepackage{amsmath,xspace,amssymb,mathrsfs}
\usepackage{color}

\input xy
\xyoption{all}
\xyoption{2cell}
\UseAllTwocells
\CompileMatrices

\newcommand{\Spec}{\operatorname{Spec}}
\renewcommand{\phi}{\varphi}

\newcommand{\Ker}{\operatorname{Ker}}

\newcommand{\Min}{\operatorname{Min}}

\newcommand{\Max}{\operatorname{Max}}

\newcommand{\Ann}{\operatorname{Ann}}

\newcommand{\Con}{\operatorname{C}}

\newtheorem{proposition}{Proposition}[section]
\newtheorem{lemma}[proposition]{Lemma} 
\newtheorem{corollary}[proposition]{Corollary}
\newtheorem{theorem}[proposition]{Theorem}

\newtheorem{prop-def}[proposition]{Proposition and definition}

\theoremstyle{definition}
\newtheorem{definition}[proposition]{Definition}
\newtheorem{example}[proposition]{Example}
\newtheorem{examples}[proposition]{Examples}
\newtheorem{remark}[proposition]{Remark}

\begin{document}

\title[Harmonic rings and lessened rings]{Structural results on harmonic rings and lessened rings}

\author[A. Tarizadeh and M. Aghajani]{Abolfazl Tarizadeh and Mohsen Aghajani}
\address{ Department of Mathematics, Faculty of Basic Sciences,
University of Maragheh \\
P. O. Box 55136-553, Maragheh, Iran.
 }
\email{ebulfez1978@gmail.com, aghajani14@gmail.com}


\begin{abstract} In this paper, a combination of algebraic and topological methods are applied to obtain new and structural results on harmonic rings. Especially, it is shown that if a Gelfand ring $A$ modulo its Jacobson radical is a zero dimensional ring, then $A$ is a clean ring. It is also proved that, for a given Gelfand ring $A$, then the retraction map $\Spec(A)\rightarrow\Max(A)$ is flat continuous if and only if $A$ modulo its Jacobson radical is a zero dimensional ring. Dually, it is proved that for a given mp-ring $A$, then the retraction map $\Spec(A)\rightarrow\Min(A)$ is Zariski continuous if and only if $\Min(A)$ is Zariski compact. New criteria for zero dimensional rings, mp-rings and Gelfand rings are given. The new notion of lessened ring is introduced and studied which generalizes ``reduced ring'' notion. Specially, a technical result is obtained which states that the product of a family of rings is a lessened ring if and only if each factor is a lessened ring. As another result in this spirit, the structure of locally lessened mp-rings is also characterized. Finally, it is characterized that a given ring when is a finite product of fields, integral domains, local rings, and lessened quasi-prime rings. \\
\end{abstract}

\maketitle

\tableofcontents

\section{Introduction and Preliminaries}

Although the present paper is the continuation of our recent work \cite{Aghajani and Tarizadeh}, but it is self-contained. For all undefined notions we refer the reader to that paper. \\

In this paper, just like that \cite{Aghajani and Tarizadeh}, a combination of algebraic and topological methods are used to solve commutative algebra problems. Specially among them, it is proved that for a given mp-ring $A$, then the retraction map $\Spec(A)\rightarrow\Min(A)$ is Zariski continuous if and only if $\Min(A)$ is Zariski compact, see Theorem \ref{Theorem Tariz-Couch}. Dually, in Theorem \ref{Theorem Tarizadeh-full}, it is shown that for a given Gelfand ring $A$, then the retraction map $\Spec(A)\rightarrow\Max(A)$ is flat continuous if and only if $A$ modulo its Jacobson radical is a zero dimensional ring. It is also proved that if a Gelfand ring $A$ modulo its Jacobson radical is a zero dimensional ring, then $A$ is a clean ring, see Theorem \ref{Theorem Harmonic VII 9294}. We introduce and study the new notion of lessened ring which generalizes reduced ring notion. Specially, a technical result is proved which states that the product of a family of rings is a lessened ring if and only if each factor is a lessened ring, see Theorem \ref{Theorem Harika}. In Theorem \ref{Lemma harmonic vii}, the structure of locally lessened mp-rings are characterized. It is also characterized that a given ring when is a finite product of fields, integral domains, local rings and lessened quasi-prime rings, see Theorems \ref{Theorem Harmonic 9462}, \ref{Corollary iv result} and \ref{Remark II} and Corollary \ref{Theorem Harmonic III}. We believe that this approach, which we call it topological commutative algebra, has the potential that can be developed further to solve more problems in commutative algebra. \\

This paper also contains new criteria for zero dimensional rings, mp-rings and Gelfand rings, see e.g. Theorems \ref{Theorem Harmonic VI 2020}, \ref{Theorem Harmonic II} and \ref{Theorem Harmonic V}. A rich literature there was dedicated to the Gelfand rings, clean rings and mp-rings (see e.g. the References and Introduction of \cite{Aghajani and Tarizadeh}). In the literature, Gelfand rings and mp-rings are also called harmonic rings (note that in the classical usage, a ring $A$ was called harmonic if $\Max(A)$ is a Hausdorff space with respect to the Zariski topology, also a Gelfand ring was called a strongly harmonic ring, see \cite[Remark 4.5]{Aghajani and Tarizadeh}, \cite[p. 99-315]{Hofmann}, \cite{Koh} and \cite{Simmons}. The word ``harmonic'' originates in Harmonic Analysis). Intuitively, the prime spectrum of a Gelfand ring can be analogized as the mountains whose the summits of the mountains are the maximal ideals, and the prime spectrum of a mp-ring can be analogized as the icicles whose the tips of the icicles are the minimal primes. \\

In this paper, all rings are commutative. If $\mathfrak{p}$ is a prime ideal of a ring $A$, then the canonical ring map $A\rightarrow A_{\mathfrak{p}}$ is denoted by $\pi_{\mathfrak{p}}$. The ideal $\Ker\pi_{\mathfrak{p}}$, especially when $\mathfrak{p}$ is a minimal prime, plays a major role in this paper. \\

A ring map is called an epimorphism of rings if it is an epimorphism in the category of commutative rings. Every surjective ring map is an epimorphism of rings, but the converse does not necessarily hold. By a flat epimorphism of rings we mean an epimorphism of rings which is also a flat ring map. \\

\section{Zero dimensionality}

\begin{lemma}\label{Lemma I Harmonic} Let $\phi: A\rightarrow B$ be a flat epimorphism of rings. If $\mathfrak{q}$ is a prime ideal of $B$, then the induced map $\phi_{\mathfrak{q}}: A_{\mathfrak{p}}\rightarrow B_{\mathfrak{q}}$ is an isomorphism of rings where $\mathfrak{p}:=\phi^{-1}(\mathfrak{q})$.\\
\end{lemma}

{\bf Proof.} If $\phi:A\rightarrow B$ is a flat ring map, then by \cite[Theorem 7.2]{Matsumura}, $\phi_{\mathfrak{q}}: A_{\mathfrak{p}}\rightarrow B_{\mathfrak{q}}$ is a faithfully flat ring map. Moreover the following diagram is commutative: $$\xymatrix{
A\ar[r]^{\phi=epic} \ar[d]^{}&B\ar[d]^{epic} \\ A_{\mathfrak{p}}\ar[r]^{\phi_{\mathfrak{q}}}&B_{\mathfrak{q}}}$$ where the vertical arrows are the canonical ring maps.
Therefore the composition $\xymatrix{A\ar[r]&A_{\mathfrak{p}} \ar[r]^{\phi_{\mathfrak{q}}}&B_{\mathfrak{q}}}$ and so $\phi_{\mathfrak{q}}$ are epimorphisms. Hence, $\phi_{\mathfrak{q}}$ is an isomorphism, because it is well known that every faithfully flat epimorphism of rings is an isomorphism. $\Box$ \\

Recall that a ring $A$ is called a $\pi-$regular ring if for every $f\in A$ there exists a natural number $n\geqslant1$ such that $f^{n}(1-fg)=0$ for some $g\in A$. This generalizes the absolutely flat ring notion. In fact, a ring is absolutely flat iff it is reduced and $\pi-$regular, because if $f^{n}(1-fg)=0$ then $f(1-fg)$ is nilpotent. \\

In the following result, new criteria for zero dimensional rings are given. Also a new proof for the implication (i)$\Rightarrow$(ii) is provided which has geometric nature, for a purely algebraic proof of this implication see \cite[Theorem 3.3]{Aghajani and Tarizadeh}. \\

\begin{theorem}\label{Theorem Harmonic VI 2020} For a ring $A$ the following statements are equivalent. \\
$\mathbf{(i)}$ $\dim(A)=0$. \\
$\mathbf{(ii)}$ Every flat epimorphism of rings
 with source $A$ is surjective. \\
$\mathbf{(iii)}$ If $\mathfrak{p}$ is a minimal prime ideal of $A$, then the canonical map $A\rightarrow A_{\mathfrak{p}}$ is surjective. \\
$\mathbf{(iv)}$ Every maximal ideal of $A$ is the radical of a unique pure ideal of $A$. \\
$\mathbf{(v)}$ $A$ is $\pi-$regular. \\
\end{theorem}

{\bf Proof.} $\mathbf{(i)}\Rightarrow\mathbf{(ii)}:$ Let $\phi:A\rightarrow B$ be a flat epimorphism of rings and let $\theta:\Spec(B)\rightarrow\Spec(A)$ be the induced morphism between the corresponding affine schemes. It suffices to show that $\theta$ is a closed immersion of schemes, because it is well known that a morphism of rings $A\rightarrow B$ is surjective if and only if the induced morphism $\Spec(B)\rightarrow\Spec(A)$ is a closed immersion of schemes. The map $\theta$ between the underlying spaces is injective, since $\phi$ is an epimorphism of rings. It is also a closed map, since $\Spec(B)$ is quasi-compact and by the hypothesis, $\Spec(A)$ is Hausdorff. It remains to show that if $\mathfrak{q}$ is a prime ideal of $B$, then $\theta_{\mathfrak{q}}^{\sharp}:\mathscr{O}_{\Spec(A),\mathfrak{p}}
\rightarrow\mathscr{O}_{\Spec(B),\mathfrak{q}}$ is surjective where $\mathfrak{p}=\theta(\mathfrak{q})=\phi^{-1}(\mathfrak{q})$. We have the following commutative diagram: $$\xymatrix{
\mathscr{O}_{\Spec(A),\mathfrak{p}}
\ar[r]^{\theta_{\mathfrak{q}}^{\sharp}}
\ar[d]&\mathscr{O}_{\Spec(B),\mathfrak{q}}
\ar[d]\\A_{\mathfrak{p}}\ar[r]^{\phi_{\mathfrak{q}}}&
B_{\mathfrak{q}}}$$ where the vertical arrows are the canonical isomorphisms. By Lemma \ref{Lemma I Harmonic}, $\phi_{\mathfrak{q}}$ is an isomorphism of rings. Therefore $\theta_{\mathfrak{q}}^{\sharp}$ is an isomorphism.
$\mathbf{(ii)}\Rightarrow\mathbf{(iii)}:$ There is nothing to prove. \\
$\mathbf{(iii)}\Rightarrow\mathbf{(i)}:$ If $\mathfrak{p}$ is a minimal prime ideal of $A$, then by the hypothesis, $A/I\simeq A_{\mathfrak{p}}$ where $I:=\Ker\pi_{\mathfrak{p}}$. So $\mathfrak{p}/I$ is the only prime ideal of $A/I$. Hence, $\mathfrak{p}$ is a maximal ideal of $A$. \\
$\mathbf{(i)}\Rightarrow\mathbf{(iv)}:$ By the hypothesis, $A$ is a mp-ring. Then apply \cite[Theorem 6.2 (x)]{Aghajani and Tarizadeh}. \\
$\mathbf{(iv)}\Rightarrow\mathbf{(i)}:$ Let $\mathfrak{p}$ be a prime ideal of $A$. There exists a maximal ideal $\mathfrak{m}$ of $A$ such that $\mathfrak{p}\subseteq\mathfrak{m}$. So there exists a pure ideal $I$ of $A$ such that $\mathfrak{m}=\sqrt{I}$.
If $f\in\mathfrak{m}$, then there exists some $g\in I$ such that $f^{n}(1-g)=0$ for some $n\geqslant1$. Hence, $f\in\mathfrak{p}$. \\
$\mathbf{(i)}\Rightarrow\mathbf{(v)}:$ By \cite[Theorem 3.3]{Aghajani and Tarizadeh}, $A/\mathfrak{N}$ is absolutely flat. So if $f\in A$ then there exists some $h\in A$ such that $f(1-fh)\in\mathfrak{N}$. Thus there exists a natural number $n\geqslant1$ such that $f^{n}(1-fg)=0$ for some $g\in A$. Hence, $A$ is $\pi-$regular. \\
$\mathbf{(v)}\Rightarrow\mathbf{(i)}:$ If $\mathfrak{p}$ is a prime ideal of $A$, then there exists a maximal ideal $\mathfrak{m}$ of $A$ such that $\mathfrak{p}\subseteq\mathfrak{m}$. If $f\in\mathfrak{m}$ then there exists some $n\geqslant1$ such that $f^{n}(1-fg)=0$. It follows that $f\in\mathfrak{p}$. Thus $\mathfrak{p}=\mathfrak{m}$. $\Box$ \\

\section{Structural results on mp-rings and lessened rings}

\begin{lemma}\label{Lemma II Harmonic} If $\mathfrak{p}$ is a prime ideal of a ring $A$, then $\sqrt{\Ker\pi_{\mathfrak{p}}}\subseteq\mathfrak{p}$. The equality holds if and only if $\mathfrak{p}$ is a minimal prime ideal of $A$. \\
\end{lemma}

{\bf Proof.} Clearly $\Ker\pi_{\mathfrak{p}}\subseteq\mathfrak{p}$ and so $\sqrt{\Ker\pi_{\mathfrak{p}}}\subseteq\mathfrak{p}$. Assume the equality holds. If $\mathfrak{q}$ is a prime ideal of $A$ such that $\mathfrak{q}\subseteq\mathfrak{p}$, then
$\Ker\pi_{\mathfrak{p}}\subseteq\mathfrak{q}$ and so $\mathfrak{p}=\mathfrak{q}$. Conversely, let $\mathfrak{p}$ be a minimal prime ideal of $A$. If $f\in\mathfrak{p}$ then there exists some $g\in A\setminus\mathfrak{p}$ such that $fg$ is nilpotent, since $\mathfrak{p}A_{\mathfrak{p}}$ is the nil-radical of $A_{\mathfrak{p}}$. It follows that $f\in\sqrt{\Ker\pi_{\mathfrak{p}}}$. $\Box$ \\

By a \emph{quasi-prime ring} we mean a ring $A$ which has a unique minimal prime ideal. Obviously this generalizes ``integral domain'' notion. Clearly a ring $A$ is a quasi-prime ring if and only if its nil-radical is a prime ideal, or equivalently, $\Spec(A)$ is an irreducible space in the Zariski topology. \\

\begin{theorem}\label{Theorem Harmonic II} For a ring $A$ the following statements are equivalent. \\
$\mathbf{(i)}$ $A$ is a mp-ring. \\
$\mathbf{(ii)}$ If $\mathfrak{p}$ and $\mathfrak{q}$ are distinct minimal prime ideals of $A$, then $\Ker\pi_{\mathfrak{p}}+\Ker\pi_{\mathfrak{q}}=A$. \\
$\mathbf{(iii)}$ $A_{\mathfrak{p}}$ is a quasi-prime ring for all $\mathfrak{p}\in\Spec(A)$. \\
$\mathbf{(iv)}$ $A_{\mathfrak{m}}$ is a quasi-prime ring for all $\mathfrak{m}\in\Max(A)$. \\
$\mathbf{(v)}$ $\sqrt{\Ker\pi_{\mathfrak{p}}}$ is a prime ideal of $A$ for all $\mathfrak{p}\in\Spec(A)$. \\
$\mathbf{(vi)}$ $\sqrt{\Ker\pi_{\mathfrak{m}}}$ is a prime ideal of $A$ for all $\mathfrak{m}\in\Max(A)$. \\
$\mathbf{(vii)}$ If $\mathfrak{p}\subseteq\mathfrak{q}$ are prime ideals of $A$, then $\sqrt{\Ker\pi_{\mathfrak{p}}}=\sqrt{\Ker\pi_{\mathfrak{q}}}$. \\
\end{theorem}

{\bf Proof.} $\mathbf{(i)}\Rightarrow\mathbf{(ii)}:$ Clearly $\mathfrak{p}+\mathfrak{q}=A$. Using Lemma \ref{Lemma II Harmonic}, then we get that $\Ker\pi_{\mathfrak{p}}+\Ker\pi_{\mathfrak{q}}=A$. \\
$\mathbf{(ii)}\Rightarrow\mathbf{(i)}:$ If $\mathfrak{p}$ and $\mathfrak{q}$ are distinct minimal prime ideals of $A$, then by the hypothesis, $\mathfrak{p}+\mathfrak{q}=A$. \\
$\mathbf{(i)}\Rightarrow\mathbf{(iii)}:$ Let $\mathfrak{q}$ be the unique minimal prime ideal of $A$ which is contained in $\mathfrak{p}$. Then $\mathfrak{q}A_{\mathfrak{p}}$ is the unique minimal prime ideal of $A_{\mathfrak{p}}$. \\ $\mathbf{(iii)}\Rightarrow\mathbf{(iv)}:$ There is nothing to prove. \\
$\mathbf{(iv)}\Rightarrow\mathbf{(i)}:$ Let $\mathfrak{p}$ and $\mathfrak{q}$ be minimal prime ideals of $A$ which are contained in a maximal ideal $\mathfrak{m}$ of $A$. Let $\mathfrak{p}'A_{\mathfrak{m}}$ be the unique minimal prime ideal of  $A_{\mathfrak{m}}$ where $\mathfrak{p}'$ is a (minimal) prime ideal of $A$. Then clearly $\mathfrak{p}'A_{\mathfrak{m}}\subseteq
\mathfrak{p}A_{\mathfrak{m}}$. It follows that $\mathfrak{p}'\subseteq\mathfrak{p}$ and so $\mathfrak{p}'=\mathfrak{p}$. Similarly we get that $\mathfrak{p}'=\mathfrak{q}$. Thus $\mathfrak{p}=\mathfrak{q}$. \\
$\mathbf{(i)}\Rightarrow\mathbf{(v)}:$ If $fg\in\sqrt{\Ker\pi_{\mathfrak{p}}}$ for some $f,g\in A$, then there exist a natural number $n\geqslant1$ and some $h\in A\setminus\mathfrak{p}$ such that $f^{n}g^{n}h=0$. Then by \cite[Theorem 6.2 (ix)]{Aghajani and Tarizadeh}, $\Ann(f^{nk})+\Ann(g^{nk}h^{k})=A$ for some $k\geqslant1$. Thus we may write $1=x+y$ where $xf^{nk}=yg^{nk}h^{k}=0$. If $x\in A\setminus\mathfrak{p}$ then $f\in\sqrt{\Ker\pi_{\mathfrak{p}}}$. Otherwise, $g\in\sqrt{\Ker\pi_{\mathfrak{p}}}$. Hence, $\sqrt{\Ker\pi_{\mathfrak{p}}}$ is a prime ideal of $A$. $\mathbf{(v)}\Rightarrow\mathbf{(vi)}:$ There is nothing to prove. \\
$\mathbf{(vi)}\Rightarrow\mathbf{(i)}:$ If $\mathfrak{p}$ is a minimal prime ideal of $A$ which is contained in a maximal ideal $\mathfrak{m}$ of $A$, then $\Ker\pi_{\mathfrak{m}}\subseteq\mathfrak{p}$ and so $\sqrt{\Ker\pi_{\mathfrak{m}}}=\mathfrak{p}$. \\
$\mathbf{(v)}\Rightarrow\mathbf{(vii)}:$ Clearly $\Ker\pi_{\mathfrak{q}}\subseteq\Ker\pi_{\mathfrak{p}}$. Now by the hypothesis and using the fact that if $\sqrt{\Ker\pi_{\mathfrak{p}}}$ is a prime ideal of a ring $A$ for some $\mathfrak{p}\in\Spec(A)$ then $\sqrt{\Ker\pi_{\mathfrak{p}}}$ is a minimal prime ideal of $A$, then we get that $\sqrt{\Ker\pi_{\mathfrak{q}}}
=\sqrt{\Ker\pi_{\mathfrak{p}}}$. \\
$\mathbf{(vii)}\Rightarrow\mathbf{(i)}:$ If $\mathfrak{p}$ and $\mathfrak{q}$ are minimal prime ideals of $A$ which are contained in a maximal ideal of $A$, then by the hypothesis, $\sqrt{\Ker\pi_{\mathfrak{p}}}=\sqrt{\Ker\pi_{\mathfrak{q}}}$. Thus by Lemma \ref{Lemma II Harmonic}, $\mathfrak{p}=\mathfrak{q}$. $\Box$ \\

\begin{corollary}\label{Corollary V} Let $A$ be a ring. If $\gamma:\Spec(A)\rightarrow\Min(A)$ is a retraction map with respect to the flat topology, then  $\gamma(\mathfrak{p})=\sqrt{\Ker\pi_{\mathfrak{p}}}$. \\
\end{corollary}

{\bf Proof.} If $\mathfrak{p}$ is a prime ideal of $A$ then
$\gamma(\mathfrak{p})\subseteq\mathfrak{p}$, because there exists a minimal prime ideal $\mathfrak{q}$ of $A$ such that $\mathfrak{q}\subseteq\mathfrak{p}$, thus
$\mathfrak{q}\in\overline{\{\mathfrak{p}\}}$ and so $\mathfrak{q}=\gamma(\mathfrak{q})\in
\overline{\{\gamma(\mathfrak{p})\}}=\{\gamma(\mathfrak{p})\}$. This also shows that $A$ is a mp-ring. Thus by Theorem \ref{Theorem Harmonic II} (v), $\gamma(\mathfrak{p})=\sqrt{\Ker\pi_{\mathfrak{p}}}$ for all $\mathfrak{p}\in\Spec(A)$. $\Box$ \\

Corollary \ref{Corollary V}, in particular, tells us that flat retraction map from $\Spec(A)$ to $\Min(A)$, if it exists, then it is unique. Note that by \cite[Theorem 6.2]{Aghajani and Tarizadeh}, there exists a flat retraction map from $\Spec(A)$ to $\Min(A)$ if and only if $A$ is a mp-ring. \\

\begin{lemma}\label{Lemma Harmonic III} If $A$ is a ring, then $\bigcap\limits_{\mathfrak{m}\in\Max(A)}\Ker\pi_{\mathfrak{m}}=0$. \\
\end{lemma}

{\bf Proof.} If $f\in\bigcap\limits_{\mathfrak{m}\in\Max(A)}\Ker\pi_{\mathfrak{m}}$ then $\Ann(f)=A$, so $f=0$. $\Box$ \\

\begin{corollary}\label{Corollary Harmonic I} For a ring $A$ the following statements are equivalent. \\
$\mathbf{(i)}$ $A$ is a reduced mp-ring. \\
$\mathbf{(ii)}$ $\Ker\pi_{\mathfrak{p}}$ is a prime ideal of $A$ for all $\mathfrak{p}\in\Spec(A)$. \\
$\mathbf{(iii)}$ $\Ker\pi_{\mathfrak{m}}$ is a prime ideal of $A$ for all $\mathfrak{m}\in\Max(A)$. \\
\end{corollary}

{\bf Proof.} $\mathbf{(i)}\Rightarrow\mathbf{(ii)}:$ By Theorem \ref{Theorem Harmonic II} (v), $\sqrt{\Ker\pi_{\mathfrak{p}}}$ is a prime ideal. But $\Ker\pi_{\mathfrak{p}}=\sqrt{\Ker\pi_{\mathfrak{p}}}$, since $A$ is reduced. $\mathbf{(ii)}\Rightarrow\mathbf{(iii)}:$ There is nothing to prove. $\mathbf{(iii)}\Rightarrow\mathbf{(i)}:$ By Theorem \ref{Theorem Harmonic II}, $A$ is a mp-ring. If $f\in A$ is nilpotent then $f\in\Ker\pi_{\mathfrak{m}}$ for all $\mathfrak{m}\in\Max(A)$. Thus by Lemma \ref{Lemma Harmonic III}, $f=0$. Hence, $A$ is reduced. $\Box$ \\

\begin{lemma}\label{Lemma Harmonic VI 56} Let $(A_{i})$ be a family of rings and $A=\prod\limits_{i}A_{i}$. If $\mathfrak{q}$ is a prime ideal of $A_{k}$, then $\Ker\pi_{\mathfrak{p}}=\pi^{-1}_{k}(\Ker\pi_{\mathfrak{q}})$
where $\mathfrak{p}:=\pi^{-1}_{k}(\mathfrak{q})$ and $\pi_{k}:A\rightarrow A_{k}$ is the canonical projection. \\
\end{lemma}

{\bf Proof.} The proof is straightforward. $\Box$ \\

\begin{proposition}\label{Proposition hamoni v} Let $(A_{i})$ be a family of rings. Then $\pi_{k}^{-1}(\mathfrak{q})$ is a minimal prime of $A=\prod\limits_{i}A_{i}$ if and only if $\mathfrak{q}$ is a minimal prime of $A_{k}$. \\
\end{proposition}

{\bf Proof.} The implication ``$\Rightarrow$'' is easy. Conversely,  suppose there exists a prime ideal $\mathfrak{p}$ of $A$ such that $\mathfrak{p}\subset\pi^{-1}_{k}(\mathfrak{q})$. Thus there exists some $a=(a_{i})\in\pi^{-1}_{k}(\mathfrak{q})$ such that $a\notin\mathfrak{p}$. Clearly $a_{k}\in\mathfrak{q}$. Thus by Lemma \ref{Lemma II Harmonic}, there exist a natural number $n\geqslant1$ and some $b_{k}\in A_{k}\setminus\mathfrak{q}$ such that $(a_{k}b_{k})^{n}=0$. Then consider the sequence $b=(b_{i})\in A$ where $b_{i}=0$ for all $i\neq k$. Then clearly $(ab)^{n}=0$. It follows that $b\in\mathfrak{p}$ and so $b_{k}\in\mathfrak{q}$ which is a contradiction. $\Box$ \\

\begin{remark} In Proposition \ref{Proposition hamoni v}, if the index set is finite then the minimal primes of $A$ are precisely of the form $\pi^{-1}_{k}(\mathfrak{p})$ with $\mathfrak{p}$ is a minimal prime of $A_{k}$. But if the index set is infinite, then the structure of the minimal primes of $A$ is complicated. For instance, the cardinality of the minimal primes of $\prod\limits_{i\geqslant1}\mathbb{Z}$ is equal to $2^{\mathfrak{c}}$ where $\mathfrak{c}$ is the cardinality of the real numbers, while the cardinality of the minimal primes of the form $\pi^{-1}_{k}(0)$ is countable, for more information see \cite[Theorem 4.5]{Tarizadeh}. \\
\end{remark}

\begin{definition} By a \emph{lessened ring} we mean a ring $A$ such that \\ $\bigcap\limits_{\mathfrak{p}\in\Min(A)}\Ker\pi_{\mathfrak{p}}=0$. We shall denote this intersection by $\mathfrak{N}'(A)$ or simply by $\mathfrak{N}'$ if there is no confusion. Clearly lessened ring generalizes reduced ring notion. The ring $\mathbb{Z}/4\mathbb{Z}$ is a lessened ring which is not reduced. \\
\end{definition}

Then we obtain the following result whose proof is technical. \\

\begin{theorem}\label{Theorem Harika} Let $(A_{i})_{i\in I}$ be a family of rings. Then $A=\prod\limits_{i}A_{i}$ is a lessened ring if and only if each $A_{i}$ is a lessened ring. \\
\end{theorem}

{\bf Proof.} Suppose $A$ is a lessened ring. If $f_{k}\in\mathfrak{N}'(A_{k})$ then it will be enough to show that $f=(f_{i})\in\mathfrak{N}'(A)$ where $f_{i}=0$ for all $i\neq k$. To see the latter it suffices to show that $f\in\Ker\pi_{\mathfrak{p}}$ for all $\mathfrak{p}\in\Min(A)$. If $\mathfrak{p}$ is a minimal prime ideal of $A$, then either $e=(\delta_{i,k})_{i\in I}\in\mathfrak{p}$ or $1-e\in\mathfrak{p}$ where $\delta_{i,k}$ is the Kronecker delta. If $e\in\mathfrak{p}$ then $e\in\Ker\pi_{\mathfrak{p}}$ and so $f=fe\in\Ker\pi_{\mathfrak{p}}$. If $1-e\in\mathfrak{p}$ then first we claim that $\pi_{k}(\mathfrak{p})$ is a prime ideal of $A$. Clearly it is an ideal of $A_{k}$, since $\pi_{k}$ is surjective. If $1\in\pi_{k}(\mathfrak{p})$ then there exists some $x=(x_{i})\in\mathfrak{p}$ such that $x_{k}=1$. It follows that $e=ex\in\mathfrak{p}$ which is a contradiction. Hence, $\pi_{k}(\mathfrak{p})$ is a proper ideal of $A$. If there exist $a,b\in A_{k}$ such that $ab\in\pi_{k}(\mathfrak{p})$ then there exists some $y=(y_{i})\in\mathfrak{p}$ such $y_{k}=ab$. Then consider the sequences $z=(z_{i})$ and $t=(t_{i})$ in $A$ where $z_{i}=1$ and $t_{i}=y_{i}$ for all $i\neq k$ and $z_{k}=a$ and $t_{k}=b$. Then clearly $zt=y\in\mathfrak{p}$. Thus either $z\in\mathfrak{p}$ or $t\in\mathfrak{p}$. It follows that either $a\in\pi_{k}(\mathfrak{p})$ or $b\in\pi_{k}(\mathfrak{p})$. This establishes the claim. So there exists a minimal prime ideal $\mathfrak{q}$ of $A_{k}$ such that
$\mathfrak{q}\subseteq\pi_{k}(\mathfrak{p})$. If $g=(g_{i})\in\pi_{k}^{-1}(\mathfrak{q})$ then there exists some $h=(h_{i})\in\mathfrak{p}$ such that $g_{k}=h_{k}$. We have $g=(1-e)g+eh\in\mathfrak{p}$. Therefore $\pi_{k}^{-1}(\mathfrak{q})\subseteq\mathfrak{p}$ and so $\pi_{k}^{-1}(\mathfrak{q})=\mathfrak{p}$. Then using Lemma \ref{Lemma Harmonic VI 56}, we have
$f\in\pi^{-1}_{k}(\Ker\pi_{\mathfrak{q}})=\Ker\pi_{\mathfrak{p}}$, since
$f_{k}\in\Ker\pi_{\mathfrak{q}}$. To see the reverse implication, it suffices to show that if $r=(r_{i})\in\mathfrak{N}'(A)$ then $r_{i}\in\mathfrak{N}'(A_{i})$ for all $i$. If $\mathfrak{q}$ is a minimal prime ideal of $A_{i}$ then by Proposition \ref{Proposition hamoni v}, $\mathfrak{p}:=\pi^{-1}_{i}(\mathfrak{q})$ is a minimal prime ideal of $A$. Therefore $r\in\Ker\pi_{\mathfrak{p}}$. Thus by Lemma \ref{Lemma Harmonic VI 56}, $r_{i}\in\Ker\pi_{\mathfrak{q}}$. $\Box$ \\

Note that if a ring has a unique prime ideal, then it is a lessened ring. In particular, for any ring $A$, then $A_{\mathfrak{p}}$ is a lessened ring for all $\mathfrak{p}\in\Min(A)$. It can be seen that $\mathbb{Z}/n\mathbb{Z}$ is a lessened ring for all $n\in\mathbb{Z}$. \\

\begin{proposition}\label{Proposition vi exam} Let $A$ be a reduced ring which has a minimal prime ideal $\mathfrak{p}$ such that it is not idempotent. Then $A/\mathfrak{p}^{2}$ is not a lessened ring. \\
\end{proposition}

{\bf Proof.} Clearly $A/\mathfrak{p}^{2}$ is a quasi-prime ring with the unique minimal prime ideal $\mathfrak{q}:=\mathfrak{p}/\mathfrak{p}^{2}$. If $f\in\mathfrak{p}$ then there exists some $g\in A\setminus\mathfrak{p}$ such that $fg=0$. It follows that $f+\mathfrak{p}^{2}\in\Ker\pi_{\mathfrak{q}}$. Hence, $\Ker\pi_{\mathfrak{q}}=\mathfrak{q}\neq0$. $\Box$ \\

\begin{examples} Every reduced noetherian local ring which is not an integral domain satisfies in the hypothesis of Proposition \ref{Proposition vi exam}. In fact, every minimal prime ideal of such ring is not idempotent. As another example, if $R$ is an integral domain then the quotient ring $R[x,y]/(xy)$ is reduced with the minimal primes $(x)/(xy)$ and $(y)/(xy)$ which none of them is idempotent. Thus by Proposition \ref{Proposition vi exam}, $R[x,y]/(x^{2},xy)$ is not a lessened ring. Moreover, the ring $R[x,y]/I$ is a quasi-prime ring with the unique minimal prime ideal $\mathfrak{p}=(x)/I$ where $I=(x^{3},x^{2}y)$, and clearly $\Ker\pi_{\mathfrak{p}}$ is neither the zero ideal nor $\mathfrak{p}$, since $0\neq x^{2}+I\in\Ker\pi_{\mathfrak{p}}$ and $x+I\in\mathfrak{p}\setminus\Ker\pi_{\mathfrak{p}}$. The latter example shows that there are rings $A$ such that the inclusions $0\subset\mathfrak{N}'\subset\mathfrak{N}$ are strict where $\mathfrak{N}$ is the nil-radical of $A$. \\
\end{examples}

\begin{proposition} For any ring $A$, then $A/\mathfrak{N}'$ is a lessened ring. \\
\end{proposition}

{\bf Proof.} The minimal primes of $A/\mathfrak{N}'$ are precisely of the form $\mathfrak{q}:=\mathfrak{p}/\mathfrak{N}'$ where $\mathfrak{p}$ is a minimal prime of $A$. We have $\Ker\pi_{\mathfrak{q}}=\Ker\pi_{\mathfrak{p}}/\mathfrak{N}'$. Thus $\mathfrak{N}'(A/\mathfrak{N}')=
\bigcap\limits_{\mathfrak{p}\in\Min(A)}
(\Ker\pi_{\mathfrak{p}}/\mathfrak{N}')=
(\bigcap\limits_{\mathfrak{p}\in\Min(A)}
\Ker\pi_{\mathfrak{p}})/\mathfrak{N}'=0$. $\Box$ \\

\begin{remark} By a mini-map we mean a ring map $\phi:A\rightarrow B$ such that $\phi^{-1}(\mathfrak{p})$ is a minimal prime ideal of $A$ for all $\mathfrak{p}\in\Min(B)$, (e.g. a flat ring map or more generally any ring map with the going-down property is a mini-map). If $\phi:A\rightarrow B$ is a mini-map, then $\mathfrak{N}'(A)\subseteq\phi^{-1}\big(\mathfrak{N}'(B)\big)$ and so we get a morphism of rings $\phi_{\mathrm{ls}}:A/\mathfrak{N}'(A)\rightarrow B/\mathfrak{N}'(B)$ given by $f+\mathfrak{N}'(A)\rightsquigarrow\phi(f)+\mathfrak{N}'(B)$. Clearly mini-maps are stable under the composition. Hence, we obtain a category whose objects are the rings and whose morphism are the mini-maps, and the assignments $A\rightsquigarrow A_{\mathrm{ls}}:=A/\mathfrak{N}'(A)$ and $\phi\rightsquigarrow\phi_{\mathrm{ls}}$ form a covariant functor from this category to the category of lessened rings (whose objects are the lessened rings and whose morphisms are the
ring maps). \\
\end{remark}

Note that if $\mathfrak{p}\subseteq\mathfrak{q}$ are prime ideals of a ring $A$, then $\Ker\pi_{P}=(\Ker\pi_{\mathfrak{p}})A_{\mathfrak{q}}$ where $P:=\mathfrak{p}A_{\mathfrak{q}}$. It follows that if $\mathfrak{m}$ is a maximal ideal of a ring $A$, then $\mathfrak{N}'(A_{\mathfrak{m}})=
\bigcap\limits_{\substack{\mathfrak{p}\in\Min(A),\\
\mathfrak{p}\subseteq\mathfrak{m}}}
(\Ker\pi_{\mathfrak{p}})A_{\mathfrak{m}}$. Unlike the reduced case, it seems that lessened rings are not stable under taking arbitrary localizations. This leads us to the following definition. We call a ring $A$ a \emph{locally lessened ring} if $A_{\mathfrak{m}}$ is a lessened ring for all $\mathfrak{m}\in\Max(A)$. \\

\begin{proposition}\label{Proposition harmoni iii} If a ring $A$ is a locally lessened ring, then $A$ is a lessened ring. \\
\end{proposition}

{\bf Proof.} Take $f\in\mathfrak{N}'(A)$. It suffices to show that $\Ann(f)=A$. If not, then there exists a maximal ideal $\mathfrak{m}$ of $A$ such that $\Ann(f)\subseteq\mathfrak{m}$. If $P$ is a minimal prime ideal of $A_{{\mathfrak{m}}}$, then there exists a minimal prime ideal $\mathfrak{p}$ of $A$ such that $\mathfrak{p}\subseteq\mathfrak{m}$
and $P=\mathfrak{p}A_{\mathfrak{m}}$. We have $\Ker\pi_{P}=(\Ker\pi_{\mathfrak{p}})A_{\mathfrak{m}}$.
Therefore $f/1\in\mathfrak{N}'(A_{\mathfrak{m}})=0$, since $A_{\mathfrak{m}}$ is a lessened ring. So there exists some $g\in A\setminus\mathfrak{m}$ such that $g\in\Ann(f)$ which is a contradiction. $\Box$ \\

\begin{theorem}\label{Lemma harmonic vii} For a ring $A$ the following statements are equivalent. \\
$\mathbf{(i)}$ $A$ is a locally lessened mp-ring. \\
$\mathbf{(ii)}$ $\Ker\pi_{\mathfrak{p}}$ is a pure ideal for all $\mathfrak{p}\in\Min(A)$. \\
$\mathbf{(iii)}$ If $\mathfrak{p}\subseteq\mathfrak{q}$ are prime ideals of $A$, then $\Ker\pi_{\mathfrak{p}}=\Ker\pi_{\mathfrak{q}}$. \\
\end{theorem}

{\bf Proof.} $\mathbf{(i)}\Rightarrow\mathbf{(ii)}:$ Take $f\in I:=\Ker\pi_{\mathfrak{p}}$. By \cite[Theorem 2.1]{Aghajani and Tarizadeh}, it suffices to show that $\Ann(f)+I=A$. If not, then there exists a maximal ideal $\mathfrak{m}$ of $A$ such that $\Ann(f)+I\subseteq\mathfrak{m}$. Thus by Theorem \ref{Theorem Harmonic II} (ii), $\mathfrak{p}\subseteq\mathfrak{m}$. We have $IA_{\mathfrak{m}}=\Ker\pi_{P}=0$ where $P:=\mathfrak{p}A_{\mathfrak{m}}$. Thus there exists some $g\in A\setminus\mathfrak{m}$ such that $g\in\Ann(f)$ which is a contradiction. \\
$\mathbf{(ii)}\Rightarrow\mathbf{(iii)}:$ Without loss of generality, we may assume that $\mathfrak{p}$ is a minimal prime of $A$. Clearly  $\Ker\pi_{\mathfrak{q}}\subseteq\Ker\pi_{\mathfrak{p}}$. If $f\in\Ker\pi_{\mathfrak{p}}$ then by the hypothesis, there exists some $g\in\Ker\pi_{\mathfrak{p}}$ such that $f(1-g)=0$. It follows that $1-g\in A\setminus\mathfrak{q}$. Therefore $f\in\Ker\pi_{\mathfrak{q}}$. \\
$\mathbf{(iii)}\Rightarrow\mathbf{(i)}:$ If $\mathfrak{p}$ and $\mathfrak{q}$ are minimal prime ideals of $A$ which are contained in a maximal ideal of $A$. Then by the hypothesis, $\Ker\pi_{\mathfrak{p}}=\Ker\pi_{\mathfrak{q}}$. Thus by Lemma \ref{Lemma II Harmonic}, $\mathfrak{p}=\mathfrak{q}$. Hence, $A$ is a mp-ring. If $\mathfrak{m}$ is a maximal ideal of $A$, then there exists a unique minimal prime ideal $\mathfrak{p}$ of $A$ such that $\mathfrak{p}\subseteq\mathfrak{m}$. Setting $\mathfrak{q}:=\mathfrak{p}A_{\mathfrak{m}}$. Then using the hypothesis we have $\Ker\pi_{\mathfrak{q}}=(\Ker\pi_{\mathfrak{p}})A_{\mathfrak{m}}=
(\Ker\pi_{\mathfrak{m}})A_{\mathfrak{m}}=0$. Hence, $A_{\mathfrak{m}}$ is a lessened ring. $\Box$ \\

\begin{corollary} If $A$ is a zero dimensional ring, then $\Ker\pi_{\mathfrak{p}}$ is a pure ideal for all $\mathfrak{p}\in\Spec(A)$. \\
\end{corollary}

{\bf Proof.} Every zero dimensional ring is a locally lessened mp-ring, then apply Theorem \ref{Lemma harmonic vii}. $\Box$ \\

\begin{proposition} If $A$ is a locally lessened mp-ring, then $A_{\mathfrak{p}}$ is a lessened ring for all $\mathfrak{p}\in\Spec(A)$. \\
\end{proposition}

{\bf Proof.} There exists a maximal ideal $\mathfrak{m}$ of $A$ such that $\mathfrak{p}\subseteq\mathfrak{m}$. By the hypothesis, there exists a unique minimal prime ideal $\mathfrak{q}$ of $A$ such that $\mathfrak{q}\subseteq\mathfrak{m}$. Clearly $\mathfrak{q}\subseteq\mathfrak{p}$. It suffices to show that $\Ker\pi_{P}=0$ where $P:=\mathfrak{q}A_{\mathfrak{p}}$. We have
$\Ker\pi_{P}=(\Ker\pi_{\mathfrak{q}})A_{\mathfrak{p}}$. But the canonical ring map $A\rightarrow A_{\mathfrak{p}}$ factors as $\xymatrix{A\ar[r]& A_{\mathfrak{m}}\ar[r]& A_{\mathfrak{p}}}$. It follows that $\Ker\pi_{P}=IA_{\mathfrak{p}}$ where $I:=(\Ker\pi_{\mathfrak{q}})A_{\mathfrak{m}}$. But $I=0$, since $A_{\mathfrak{m}}$ is a lessened ring. Hence, $\Ker\pi_{P}=0$. $\Box$ \\

We call a ring $A$ a purified ring if $\mathfrak{p}$ and $\mathfrak{q}$ are distinct minimal prime ideals of $A$, then there exists an idempotent $e\in A$ such that $e\in\mathfrak{p}\setminus\mathfrak{q}$. Clearly every purified ring is a mp-ring. Purified rings are quite interesting, for more information see \cite[\S8]{Aghajani and Tarizadeh}. Specially in \cite[Theorem 8.5]{Aghajani and Tarizadeh}, reduced purified rings are characterized. One of the main motivations in defining lessened ring notion is to try to characterize general purified rings (not necessarily reduced). By a regular ideal we mean an ideal which is generated by a set of idempotents. \\

\begin{corollary}\label{Corollary harmonic II 2020} Let $A$ be a ring. If  $\Ker\pi_{\mathfrak{p}}$ is a regular ideal for all $\mathfrak{p}\in\Min(A)$, then $A$ is a locally lessened purified ring. \\
\end{corollary}

{\bf Proof.} By Theorem \ref{Lemma harmonic vii}, $A$ is a locally lessened ring, since every regular ideal is a pure ideal. If $\mathfrak{p}$ and $\mathfrak{q}$ are distinct minimal prime ideals of $A$ then by \cite[Lemma 3.2]{Aghajani and Tarizadeh},
there exists $f\in A\setminus\mathfrak{p}$ and $g\in A\setminus\mathfrak{q}$ such that $fg=0$. It follows that $g\in\Ker\pi_{\mathfrak{p}}$. So by \cite[Lemma 8.4]{Aghajani and Tarizadeh}, there exists an idempotent $e\in\Ker\pi_{\mathfrak{p}}\subseteq\mathfrak{p}$ such that $g(1-e)=0$. This yields that $1-e\in\mathfrak{q}$. Hence, $A$ is also a purified ring. $\Box$ \\

\begin{proposition} Let $A$ be a ring. Then $A$ is a lessened ring if and only if $\mathfrak{N}'$ is a pure ideal. \\
\end{proposition}

{\bf Proof.} Clearly $\sqrt{0}=\sqrt{\mathfrak{N}'}$. If $\mathfrak{N}'$ is a pure ideal then $\mathfrak{N}'=0$, because if $I$ and $J$ are pure ideals of a ring $A$ such that $\sqrt{I}=\sqrt{J}$ then $I=J$. As another reason, the only pure ideal contained in the Jacobson radical is the zero ideal.
 $\Box$ \\

\begin{lemma}\label{Lemma Harmonic IV} If $\mathfrak{p}$ is a minimal prime ideal of a ring $A$, then $A/\Ker\pi_{\mathfrak{p}}$ is a lessened quasi-prime ring. \\
\end{lemma}

{\bf Proof.} Let $\mathfrak{q}/\Ker\pi_{\mathfrak{p}}$ be a minimal prime ideal of $A/\Ker\pi_{\mathfrak{p}}$ where $\mathfrak{q}$ is a prime ideal of $A$ such that $\Ker\pi_{\mathfrak{p}}\subseteq\mathfrak{q}$. Then by Lemma \ref{Lemma II Harmonic}, $\mathfrak{p}=\sqrt{\Ker\pi_{\mathfrak{p}}}\subseteq\mathfrak{q}$ and so $\mathfrak{p}/\Ker\pi_{\mathfrak{p}}\subseteq
\mathfrak{q}/\Ker\pi_{\mathfrak{p}}$. It follows that $\mathfrak{p}=\mathfrak{q}$. As another proof, $A/\Ker\pi_{\mathfrak{p}}$ can be viewed as a subring of $A_{\mathfrak{p}}$. Thus by \cite[Lemma 3.9]{Abolfazl}, $A/\Ker\pi_{\mathfrak{p}}$ is a quasi-prime ring. Then we show that $\Ker\pi_{P}=0$ where $P=\mathfrak{p}/\Ker\pi_{\mathfrak{p}}$. If $f+\Ker\pi_{\mathfrak{p}}\in\Ker\pi_{P}$ then there exists some $g\in A\setminus\mathfrak{p}$ such that $fg\in\Ker\pi_{\mathfrak{p}}$. So $fgh=0$ for some $h\in A\setminus\mathfrak{p}$. It follows that $f\in\Ker\pi_{\mathfrak{p}}$. $\Box$ \\

\begin{lemma}\label{Lemma viii har} If $I$ is a finitely generated and idempotent ideal of a ring $A$, then there exists a unique idempotent $f\in A$ such that $I=Af$. \\
\end{lemma}

{\bf Proof.} There exists some $f\in I$ such that $(1-f)I=0$. It follows that $f$ is an idempotent and $I=Af$. The uniqueness is obvious. $\Box$ \\

\begin{theorem}\label{Theorem Harmonic 9462} For a ring $A$ the following statements are equivalent. \\
$\mathbf{(i)}$ $A$ is a finite product of lessened quasi-prime rings. \\
$\mathbf{(ii)}$ $A$ is a lessened mp-ring with finitely many minimal primes. \\
$\mathbf{(iii)}$ $\Ker\pi_{\mathfrak{p}}$ is a finitely generated and idempotent ideal of $A$ for all $\mathfrak{p}\in\Min(A)$. \\
\end{theorem}

{\bf Proof.} $\mathbf{(i)}\Rightarrow\mathbf{(ii)}:$ Clearly $A$ is a mp-ring, because the product of a finite family of rings is a mp-ring if and only if each factor is a mp-ring. By the hypothesis, $A=A_{1}\times...\times A_{n}$ where each $A_{i}$ is a lessened and quasi-prime ring. If $\mathfrak{p}$ is a minimal prime of $A$, then $\mathfrak{p}=A_{1}\times...\times\mathfrak{p}_{k}\times...\times A_{n}$ where $\mathfrak{p}_{k}$ is a minimal prime of $A_{k}$ for some $k$. Thus by Lemma \ref{Lemma Harmonic VI 56},
$\Ker\pi_{\mathfrak{p}}=A_{1}
\times...\times I_{k}\times...\times A_{n}$ where $I_{k}=\{0\}$.
Therefore $\mathfrak{N}'(A)
=\bigcap\limits_{k=1}^{n}(A_{1}
\times...\times I_{k}\times...\times A_{n})=0$. \\
$\mathbf{(ii)}\Rightarrow\mathbf{(i)}:$ If $\mathfrak{p}$ and $\mathfrak{q}$ are distinct minimal primes of $A$, then by Theorem \ref{Theorem Harmonic II} (ii), $\Ker\pi_{\mathfrak{p}}+\Ker\pi_{\mathfrak{q}}=A$. Thus by the Chinese remainder theorem,  $A\simeq A/\Ker\pi_{\mathfrak{p}_{1}}\times...\times A/\Ker\pi_{\mathfrak{p}_{n}}$ where $\Min(A)=\{\mathfrak{p}_{1},...,\mathfrak{p}_{n}\}$. Then by applying Lemma \ref{Lemma Harmonic IV}, the assertion is concluded. \\
$\mathbf{(i)}\Rightarrow\mathbf{(iii)}:$ By the hypothesis and Lemma \ref{Lemma Harmonic VI 56}, there exists some $k$ such that $\Ker\pi_{\mathfrak{p}}$ is generated by the idempotent $1_{A}-(\delta_{i,k})_{i=1}^{n}=(1,...,0,...,1)$. Hence, $\Ker\pi_{\mathfrak{p}}$ is a finitely generated and idempotent ideal. \\
$\mathbf{(iii)}\Rightarrow\mathbf{(ii)}:$ If $\mathfrak{p}\in\Min(A)$ then by Lemma \ref{Lemma viii har}, there exists an idempotent  $e_{\mathfrak{p}}\in A$ such that $\Ker\pi_{\mathfrak{p}}=Ae_{\mathfrak{p}}$.
Thus by Theorem \ref{Lemma harmonic vii}, $A$ is a locally lessened mp-ring. So by Proposition \ref{Proposition harmoni iii}, $A$ is a lessened ring. We also have $\Spec(A)=\bigcup\limits_{\mathfrak{p}\in\Min(A)}
D(1-e_{\mathfrak{p}})$. Using the quasi-compactness of the prime spectrum, then we may write $\Spec(A)=\bigcup\limits_{i=1}^{n}
D(1-e_{i})$ where $e_{i}:=e_{\mathfrak{p}_{i}}$ for all $i$.
If $\mathfrak{q}$ is a minimal prime of $A$, then $e_{\mathfrak{p}}\in\mathfrak{q}$ for some $\mathfrak{p}\in\{\mathfrak{p}_{1},...,\mathfrak{p}_{n}\}$. It follows that $\Ker\pi_{\mathfrak{p}}=Ae_{\mathfrak{p}}\subseteq\mathfrak{q}$. Thus by Lemma \ref{Lemma II Harmonic}, $\mathfrak{p}=\sqrt{\Ker\pi_{\mathfrak{p}}}
\subseteq\mathfrak{q}$. So $\mathfrak{p}=\mathfrak{q}$. Hence, $\Min(A)=\{\mathfrak{p}_{1},...,\mathfrak{p}_{n}\}$.
$\Box$ \\

\begin{theorem}\label{Corollary iv result} For a ring $A$ the following statements are equivalent. \\
$\mathbf{(i)}$ $A$ is a finite product of integral domains. \\
$\mathbf{(ii)}$ $A$ is a reduced mp-ring with finitely many minimal primes. \\
$\mathbf{(iii)}$ Every minimal prime of $A$ is a finitely generated and idempotent ideal. \\
\end{theorem}

{\bf Proof.} It is proved exactly like Theorem \ref{Theorem Harmonic 9462}. $\Box$ \\

\begin{corollary}\label{Theorem Harmonic III} Let $A$ be a ring. Then $A$ is a finite product of fields if and only if every maximal ideal of $A$ is a finitely generated and idempotent ideal. \\
\end{corollary}

{\bf Proof.} The implication ``$\Rightarrow$" is easy. Conversely, if $\mathfrak{m}$ is a maximal ideal of $A$ then there exists an idempotent $e\in A$ such that $\mathfrak{m}=Ae$. Hence, $A$ is a zero dimensional ring. Thus by Theorem \ref{Corollary iv result} (iii), $A$ is a finite product of integral domains. Each of these integral domains is a field, since $A$ is a zero dimensional ring. $\Box$ \\

In the following result, a new proof is given to \cite[Theorem 6.3]{Aghajani and Tarizadeh}. \\

\begin{theorem}\label{Theorem Tariz-Couch} Let $A$ be a mp-ring. Then the retraction map $\gamma:\Spec(A)\rightarrow\Min(A)$ is Zariski continuous if and only if $\Min(A)$ is Zariski compact. \\
\end{theorem}

{\bf Proof.} The implication ``$\Rightarrow$'' is easy. To see the reverse implication, it suffices to show that $\gamma^{-1}\big(W_{f})$ is a Zariski open of $\Spec(A)$ for all $f\in A$, where $W_{f}:=\Min(A)\cap D(f)$.
It can be shown that $W_{f}$ is a Zariski clopen of $\Min(A)$.
By \cite[Theorem 4.3]{Tarizadeh}, the induced Zariski and flat topologies over $\Min(A)$ are the same. Hence, $W_{f}$ is a flat clopen of $\Min(A)$. So $\gamma^{-1}\big(W_{f})$ is a flat clopen of $\Spec(A)$, because by \cite[Theorem 6.2 (v)]{Aghajani and Tarizadeh}, the map $\gamma$ is flat continuous. By \cite[Corollary 3.12]{Abolfazl}, the Zariski clopens and flat clopens of $\Spec(A)$ are the same. $\Box$ \\

\begin{proposition}\label{new Prop 1} Let $A$ be a ring. Then $A$ is a mp-ring if and only if the polynomial ring $A[x]$ is a mp-ring. \\
\end{proposition}

{\bf Proof.} The assertion is deduced from the general fact that for any ring $A$, then the minimal prime ideals of $A[x]$ are precisely of the form $\mathfrak{p}[x]$ where $\mathfrak{p}$ is a minimal prime ideal of $A$. $\Box$ \\

In the same vein, see also Proposition \ref{Proposition Harmonic I}. \\

\section{New results on Gelfand rings and clean rings}

In the following results, $\mathfrak{J}$ denotes the Jacobson radical of $A$. \\

\begin{lemma}\label{Lemma Harmonic V 20} For a ring $A$ the following statements are equivalent. \\
$\mathbf{(i)}$ $A/\mathfrak{J}$ is a zero dimensional ring. \\
$\mathbf{(ii)}$ $\Max(A)\cap V(f)$ is a Zariski open of $\Max(A)$ for all $f\in A$. \\
\end{lemma}

{\bf Proof.} $\mathbf{(i)}\Rightarrow\mathbf{(ii)}:$ There exists some $g\in A$ such that $f(1-fg)\in\mathfrak{J}$, because $A/\mathfrak{J}$ is reduced and so it is an absolutely flat ring. It follows that $\Max(A)\cap V(f)=\Max(A)\cap D(1-fg)$.
$\mathbf{(ii)}\Rightarrow\mathbf{(i)}:$ Let $\mathfrak{p}$ be a prime ideal of $A$ such that $\mathfrak{J}\subseteq\mathfrak{p}$. There exists a maximal ideal $\mathfrak{m}$ of $A$ such that $\mathfrak{p}\subseteq\mathfrak{m}$. If $f\in\mathfrak{m}$ then by the hypothesis, there exists some $g\in A$ such that $\mathfrak{m}\in\Max(A)\cap D(g)\subseteq\Max(A)\cap V(f)$. It follows that $fg\in\mathfrak{J}$ and so $f\in\mathfrak{p}$. Hence, $\mathfrak{p}=\mathfrak{m}$. $\Box$ \\

Then we obtain the following non-trivial result. \\

\begin{theorem}\label{Theorem Harmonic VII 9294} Let $A$ be a Gelfand ring such that $A/\mathfrak{J}$ is zero dimensional ring. Then $A$ is a clean ring. \\
\end{theorem}

{\bf Proof.} By Lemma \ref{Lemma Harmonic V 20}, $\Max(A)\cap V(f)$ is a Zariski clopen of $\Max(A)$ for all $f\in A$. It follows that $\Max(A)$ is totally disconnected with respect to the Zariski topology. Therefore by \cite[Theorem 5.1 (vi)]{Aghajani and Tarizadeh}, $A$ is a clean ring. $\Box$ \\

\begin{theorem}\label{Theorem Tarizadeh-full} Let $A$ be a Gelfand ring. Then the retraction map $\gamma:\Spec(A)\rightarrow\Max(A)$ is flat continuous if and only if $A/\mathfrak{J}$ is a zero dimensional ring. \\
\end{theorem}

{\bf Proof.} Assume $\gamma$ is flat continuous. If $f\in A$ then $\gamma^{-1}(U_{f})$ is a flat clopen of $\Spec(A)$, because for any ring $A$, then $U_{f}:=\Max(A)\cap V(f)$ is a flat clopen of $\Max(A)$. But the Zariski clopens and flat clopens of $\Spec(A)$ are the same, see \cite[Corollary 3.12]{Abolfazl}. Thus there exists an idempotent $e\in A$ such that $\gamma^{-1}(U_{f})=D(e)$. It follows that $U_{f}=\Max(A)\cap D(e)$. Thus by Lemma \ref{Lemma Harmonic V 20}, $A/\mathfrak{J}$ is a zero dimensional ring. To see the reverse implication, it suffices to show that $\gamma^{-1}(U_{f})$ is a flat open of $\Spec(A)$, because the collection of $U_{f}$ with $f\in A$ forms a subbase for the induced flat topology over $\Max(A)$. By Lemma \ref{Lemma Harmonic V 20}, $U_{f}$ is a Zariski clopen of $\Max(A)$. By \cite[Theorem 4.3 (ii)]{Aghajani and Tarizadeh}, the map $\gamma$ is Zariski continuous, thus there exists an idempotent $e\in A$ such that $\gamma^{-1}(U_{f})=V(e)$. $\Box$ \\

\begin{theorem}\label{Remark II} For a ring $A$ the following statements are equivalent. \\
$\mathbf{(i)}$ $A$ is a finite product of local rings. \\
$\mathbf{(ii)}$ $A$ is a Gelfand ring with finitely many maximal ideals. \\
$\mathbf{(iii)}$ $A$ is a Gelfand ring and $\Ker\pi_{\mathfrak{m}}$ is a finitely generated and idempotent ideal of $A$ for all $\mathfrak{m}\in\Max(A)$. \\
\end{theorem}

{\bf Proof.} $\mathbf{(i)}\Rightarrow\mathbf{(ii)}:$ By \cite[Theorem 4.3 (iv)]{Aghajani and Tarizadeh}, the product of a family of rings is a Gelfand ring if and only if each factor is a Gelfand ring. \\
 $\mathbf{(ii)}\Rightarrow\mathbf{(i)}:$ If $\mathfrak{m}$ and $\mathfrak{m}'$ are distinct maximal ideals of $A$ then by \cite[Theorem 4.3 (ix)]{Aghajani and Tarizadeh}, $\Ker\pi_{\mathfrak{m}}+\Ker\pi_{\mathfrak{m}'}=A$.
Then using Lemma \ref{Lemma Harmonic III} and the Chinese remainder theorem, we have $A\simeq A/\Ker\pi_{\mathfrak{m}_{1}}\times...\times
A/\Ker\pi_{\mathfrak{m}_{d}}$ where $\Max(A)=\{\mathfrak{m}_{1},...,\mathfrak{m}_{d}\}$. \\
$\mathbf{(i)}\Rightarrow\mathbf{(iii)}:$ By the hypothesis, $A=A_{1}\times...\times A_{d}$ where each $A_{i}$ is a local ring with the maximal ideal $\mathfrak{m}_{i}$. If $\mathfrak{m}$ is a maximal ideal of $A$ then $\mathfrak{m}=A_{1}\times...\times\mathfrak{m}_{k}\times... \times A_{d}$ for some $k$. It follows that $\Ker\pi_{\mathfrak{m}}$ is generated by the idempotent $1_{A}-(\delta_{i,k})_{i=1}^{d}=(1,...,0,...,1)$. Hence,  $\Ker\pi_{\mathfrak{m}}$ is a finitely generated and idempotent ideal. \\
$\mathbf{(iii)}\Rightarrow\mathbf{(ii)}:$ By Lemma \ref{Lemma viii har}, there exists an idempotent $e_{\mathfrak{m}}\in A$ such that $\Ker\pi_{\mathfrak{m}}=Ae_{\mathfrak{m}}$. It follows that $\Spec(A)=\bigcup\limits_{\mathfrak{m}\in\Max(A)}
D(1-e_{\mathfrak{m}})$. Using the quasi-compactness of the prime spectrum, then we may write $\Spec(A)=\bigcup\limits_{i=1}^{d}
D(1-e_{i})$ where $e_{i}:=e_{\mathfrak{m}_{i}}$ for all $i$.
If $\mathfrak{m}$ is a maximal ideal of $A$, then $e_{\mathfrak{m}'}\in\mathfrak{m}$ for some $\mathfrak{m}'\in\{\mathfrak{m}_{1},...,\mathfrak{m}_{d}\}$. It follows that $\Ker\pi_{\mathfrak{m}'}=Ae_{\mathfrak{m}'}\subseteq\mathfrak{m}$. Thus by \cite[Theorem 4.3 (ix)]{Aghajani and Tarizadeh}, $\mathfrak{m}=\mathfrak{m}'$. Hence, $\Max(A)=\{\mathfrak{m}_{1},...,\mathfrak{m}_{d}\}$. $\Box$ \\

\begin{proposition}\label{Proposition Harmonic I} Let $A$ be a ring. If $A[x]$ is a Gelfand ring, then $A$ is a Gelfand ring. \\
\end{proposition}

{\bf Proof.} If a prime ideal $\mathfrak{p}$ of $A$ is contained in distinct maximal ideals $\mathfrak{m}$ and $\mathfrak{m}'$ of $A$, then there exists a maximal ideal $M$ of $A[x]$ such that $\mathfrak{m}[x]+\mathfrak{m}'[x]\subseteq M$. But this is a contradiction, since $\mathfrak{m}+\mathfrak{m}'=A$ and so $\mathfrak{m}[x]+\mathfrak{m}'[x]=A[x]$. $\Box$ \\

The converse of the above result is not true. For example, if $k$ is a field then $k[x]$ is not a Gelfand ring. If $A[x]$ is a clean ring, then $A$ is a clean ring. \\

\begin{proposition}\label{Corollary pure not regular} If $X$ is a connected topological space with at least two distinct points then $\Con(X)$, the ring of real-valued continuous functions on $X$, is a Gelfand ring which is not a clean ring. \\
\end{proposition}

{\bf Proof.} It is well known that for any topological space $X$, then $\Con(X)$ is a Gelfand ring. Also for each $x\in X$, $\mathfrak{m}_{x}=\{f\in\Con(X): f(x)=0\}$ is a maximal ideal of $\Con(X)$. It is easy to see that a space $X$ is connected iff $\Con(X)$ has no nontrivial idempotents. Thus by \cite[Theorem 5.1 (ix)]{Aghajani and Tarizadeh}, $\Con(X)$ is not a clean ring. $\Box$ \\

\begin{example} (A pure ideal which is not regular) Finding a pure ideal which is not a regular ideal, is not as easy as one may think at first. Proposition \ref{Corollary pure not regular} together with  \cite[Theorem 5.1(iii)]{Aghajani and Tarizadeh} guarantee that the existence of a pure ideal in $\Con(X)$ which is not a regular ideal. \\
\end{example}

If $A$ is a Gelfand ring, then by \cite[Theorem 4.3 (vi)]{Aghajani and Tarizadeh},  $\Ker\pi_{\mathfrak{m}}$ is a pure ideal for all $\mathfrak{m}\in\Max(A)$. But the converse does not necessarily hold. For example, take the ring of integers. \\

If $I$ is an ideal of a ring $A$, then we define $I^{\ast}=\{f\in A: Af+I\neq A\}$. If $I$ is a proper ideal of $A$, then $I\subseteq I^{\ast}$. Clearly $$I^{\ast}=\bigcup\limits_{\mathfrak{p}\in V(I)}\mathfrak{p}=\bigcup\limits_{\mathfrak{m}\in V(I)\cap\Max(A)}\mathfrak{m}.$$ In general, $I^{\ast}$ is not an ideal of $A$. For example, if $I$ is the zero ideal of the ring of integers, then $-2,3\in I^{\ast}=\{f\in\mathbb{Z}:f\mathbb{Z}
\neq\mathbb{Z}\}$ but $-2+3=1\notin I^{\ast}$. \\

\begin{lemma} Let $I$ be an ideal of a ring $A$. Then $I=I^{\ast}$ if and only if $I$ is a maximal ideal of $A$. \\
\end{lemma}

{\bf Proof.} If $I=I^{\ast}$ then $I$ is a proper ideal of $A$. So there exists a maximal ideal $\mathfrak{m}$ of $A$ such that $I\subseteq\mathfrak{m}$. It follows that $\mathfrak{m}\subseteq\mathfrak{m}^{\ast}\subseteq I^{\ast}$. Thus $I=\mathfrak{m}$. The reverse implication is easily proved. $\Box$ \\

\begin{proposition}\label{Proposition Harmonic II} Let $I$ be an ideal of a ring $A$. Then $I^{\ast}$ is an ideal of $A$ if and only if $A/I$ is a local ring. \\
\end{proposition}

{\bf Proof.} If $I^{\ast}$ is an ideal of $A$ then $I$ is a proper ideal of $A$, since $0\in I^{\ast}$. Hence, $A/I$ is a non-zero ring. Suppose $I$ is contained in distinct maximal ideals $\mathfrak{m}$ and $\mathfrak{m}'$ of $A$.  Then we have $\mathfrak{m}=\mathfrak{m}^{\ast}\subseteq I^{\ast}$. Similarly we get that $\mathfrak{m}'\subseteq I^{\ast}$. Thus $A=\mathfrak{m}+\mathfrak{m}'\subseteq I^{\ast}$, since $I^{\ast}$ is an ideal of $A$. So $1\in I^{\ast}$ which is a contradiction. Conversely, let $\mathfrak{m}$ be the unique maximal ideal of $A$ containing $I$. Then $I^{\ast}=\mathfrak{m}$. $\Box$ \\

\begin{corollary} For a non-zero ring $A$ the following statements are equivalent. \\
$\mathbf{(i)}$ $A$ is a local ring. \\
$\mathbf{(ii)}$ If $I$ is a proper ideal of $A$, then $I^{\ast}$ is an ideal of $A$. \\
$\mathbf{(iii)}$ $\{0\}^{\ast}$ is an ideal of $A$. \\
\end{corollary}

{\bf Proof.} It follows from Proposition \ref{Proposition Harmonic II}. $\Box$ \\

\begin{theorem}\label{Theorem Harmonic V} For a ring $A$ the following statements are equivalent. \\
$\mathbf{(i)}$ $A$ is a Gelfand ring. \\
$\mathbf{(ii)}$ $A/\Ker\pi_{\mathfrak{p}}$ is a local ring for all $\mathfrak{p}\in\Spec(A)$. \\
$\mathbf{(iii)}$ $\mathfrak{p}^{\ast}$ is an ideal of $A$ for all $\mathfrak{p}\in\Spec(A)$. \\
$\mathbf{(iv)}$ $\mathfrak{p}^{\ast}$ is an ideal of $A$ for all $\mathfrak{p}\in\Min(A)$. \\
\end{theorem}

{\bf Proof.} $\mathbf{(i)}\Rightarrow\mathbf{(ii)}:$ There exists a maximal ideal $\mathfrak{m}$ of $A$ such that $\mathfrak{p}\subseteq\mathfrak{m}$. Thus $\Ker\pi_{\mathfrak{m}}\subseteq\Ker\pi_{\mathfrak{p}}$. By
\cite[Theorem 4.3 (ix)]{Aghajani and Tarizadeh}, $A/\Ker\pi_{\mathfrak{m}}$ is a local ring. Hence, $A/\Ker\pi_{\mathfrak{p}}$ is a local ring.
$\mathbf{(ii)}\Rightarrow\mathbf{(i)}:$ Easy.
$\mathbf{(i)}\Rightarrow\mathbf{(iii)}:$ By the hypothesis, $A/\mathfrak{p}$ is a local ring and so by Proposition \ref{Proposition Harmonic II}, $\mathfrak{p}^{\ast}$ is an ideal of $A$.  $\mathbf{(iii)}\Rightarrow\mathbf{(iv)}:$ There is nothing to prove.
$\mathbf{(iv)}\Rightarrow\mathbf{(i)}:$ If $\mathfrak{p}$ is a minimal prime ideal of $A$, then by Proposition \ref{Proposition Harmonic II}, $A/\mathfrak{p}$ is a local ring. $\Box$ \\

The following result is the dual of Corollary \ref{Corollary V}. \\

\begin{corollary} Let $A$ be a ring. If $\gamma:\Spec(A)\rightarrow\Max(A)$ is a retraction map with respect to the Zariski topology, then  $\gamma(\mathfrak{p})=\mathfrak{p}^{\ast}$. \\
\end{corollary}

{\bf Proof.} If $\mathfrak{p}$ is a prime ideal of $A$ then $\mathfrak{p}\subseteq\gamma(\mathfrak{p})$, because there exists a maximal ideal $\mathfrak{m}$ of $A$ such that $\mathfrak{p}\subseteq\mathfrak{m}$, thus $\mathfrak{m}\in V(\mathfrak{p})=\overline{\{\mathfrak{p}\}}$ and so $\mathfrak{m}=\gamma(\mathfrak{m})\in
\overline{\{\gamma(\mathfrak{p})\}}=\{\gamma(\mathfrak{p})\}$, hence $\mathfrak{m=\gamma(\mathfrak{p})}$. Therefore $A$ is a Gelfand ring and so by Theorem \ref{Theorem Harmonic V}, $\gamma(\mathfrak{p})=\mathfrak{p}^{\ast}$. $\Box$ \\

\section{Some dual aspects of harmonic rings}

Let $A$ be a ring. By \cite[Theorem 6]{Hochster} or \cite[Theorem 3.20]{Abolfazl}, there exists a ring $B$ and a homeomorphism $\phi:\big(\Spec(A),\mathcal{Z}\big)\rightarrow
\big(\Spec(B),\mathcal{F}\big)$ such that if $\mathfrak{p}\subseteq\mathfrak{q}$ are primes of $A$, then $\phi(\mathfrak{q})\subseteq\phi(\mathfrak{p})$ where $\mathcal{Z}$ (resp. $\mathcal{F}$) denotes the Zariski (resp. flat) topology. We call such ring $B$ a \emph{prime-inverse ring} of $A$. In this case, it can be seen that the map $\phi^{-1}:\big(\Spec(B),\mathcal{Z}\big)\rightarrow
\big(\Spec(A),\mathcal{F}\big)$ is also a homeomorphism. Hence, to have such homeomorphisms, does not matter $\Spec(A)$ or $\Spec(B)$ with which topologies (Zariski or flat) are equipped, it suffices only to consider the mutual topologies over them. Therefore the prime-inverse notion is symmetric. Note that prime-inverse ring is not unique. Clearly every prime-inverse of a Gelfand ring is a mp-ring, (and vise versa). \\

\begin{theorem}\label{Theorem Harmonic I} Every prime-inverse of a clean ring is a purified ring. \\
\end{theorem}

{\bf Proof.} Let $A$ be a clean ring and let $B$ be a prime-inverse ring of $A$. Let $\mathfrak{p}$ and $\mathfrak{p}'$ be distinct minimal prime ideals of $B$. There exists a homeomorphism $\phi:\Spec(A)\rightarrow\Spec(B)$ which reverses inclusions. Thus there exist maximal ideals $\mathfrak{m}$ and $\mathfrak{m}'$ of $A$ such that $\phi(\mathfrak{m})=\mathfrak{p}$ and $\phi(\mathfrak{m}')=\mathfrak{p}'$. By \cite[Theorem 5.1 (ix)]{Aghajani and Tarizadeh}, there exists an idempotent $e\in A$ such that $e\in\mathfrak{m}$ and $1-e\in\mathfrak{m}'$. Then by \cite[Theorem 2.2]{Aghajani and Tarizadeh}, there exists an idempotent $e'\in B$ such that $\phi\big(D(e)\big)=D(e')$. It follows that $e'\in\mathfrak{p}$ and $1-e'\in\mathfrak{p}'$. $\Box$ \\

\begin{theorem} Every prime-inverse of a purified ring is a clean ring. \\
\end{theorem}

{\bf Proof.} It is proved exactly like Theorem \ref{Theorem Harmonic I}. $\Box$ \\

\end{document}